\documentclass[reqno]{amsart}
\usepackage{amsmath, amsfonts, amsthm, amssymb}
\usepackage{mathrsfs}
\usepackage[a4paper]{geometry}
\usepackage{graphicx}
\newtheorem{thm}{Theorem}[section]
\newtheorem{prop}[thm]{Proposition}
\newtheorem{lem}[thm]{Lemma}
\newtheorem{cor}[thm]{Corollary}
\theoremstyle{definition}
\newtheorem{df}[thm]{Definition}
\theoremstyle{remark}
\newtheorem{expl}{Example}[section]
\newtheorem*{explH}{Example (heat wavelet families)}
\newtheorem{rem}[expl]{Remark}
\numberwithin{equation}{section}
\def\R{\mathbb R}
\def\Z{\mathbb Z}
\def\N{\mathbb N}
\def\C{\mathbb C}
\def\S{\mathbb S}
\def\I{\mathrm I}
\def\d{\mathrm d}
\def\G{\mathscr G}
\def\H{\mathscr H}
\def\e{\mathrm e}
\def\i{\mathrm i}
\DeclareMathOperator{\trace}{trace}
\DeclareMathOperator{\spann}{span}
\DeclareMathOperator{\bdiag}{b-diag}
\DeclareMathOperator{\rank}{rank}
\def\SO{\mathrm{SO}}
\def\U{\mathrm U}

\DeclareMathOperator{\diag}{diag}
\title{Diffusive wavelets on groups and homogeneous spaces}
\author{Svend Ebert}
\author{Jens Wirth}\thanks{J.Wirth is supported by the Engineering and Physics Research Council, UK, by grant
EP/E062873/1.} 
\begin{document}

\begin{abstract}
The aim of this exposition is to explain basic ideas behind the concept of diffusive wavelets on spheres in the language of representation theory of Lie groups and within the framework of the group Fourier transform given by Peter-Weyl decomposition of $L^2(\G)$ for a compact Lie group $\G$. 

After developing a general concept for compact groups and their homogeneous spaces we give concrete examples for tori---which reflect the situation on $\R^n$---and for spheres $\S^2$ and $\S^3$.
\end{abstract}

\maketitle

\tableofcontents

\section{Introduction}

Within the past decades wavelets and associated wavelet transforms have been intensively used in both applied and pure mathematics. They and the related multi-scale analysis provide essential tools to describe, analyse and modify signals, images or in rather abstract concepts functions, function spaces and associated operators. 
We refer to the seminal expositions of Meyer \cite{Meyer} and Coifman--Meyer \cite{CoifmanMeyer}. The most influential constructions in $\R^n$ are due to Haar \cite{Haar}, Grossmann--Morlet \cite{GrossmannMorlet} and Daubechies \cite{Daubechies}. 

The constructions of wavelet transforms can be entirely based on an abstract group theoretic and representation theoretic approach. An abstract exposition of this can be found in Kisil \cite{Kisil}, for the particular situation of the Lorentz group $\SO(n+1,1)$ acting on the sphere $\S^n$ an associated associated wavelet constructions was carried out by Antoine and Vandergheynst \cite{AntoineVandergheynst}, \cite{AntoineVandergheynst2}; see also \cite{Ferreira}, \cite{Ferreira2}.

We will shortly recall the general group theoretic approach. Let $\mathcal M$ be a (pseudo-) Riemannian manifold. A wavelet transform in $L^2(\mathcal M)$ is defined in terms of an irreducible unitary representation $U$ of Lie group $\mathscr G$ 
\begin{equation}
   U : \mathscr G \to \mathcal L( L^2(\mathcal M)).
\end{equation}
A non-zero vector $\Psi\in L^2(\mathcal M)$ is an admissible wavelet if
\begin{equation}\label{eq:admCond}
  \int_{\mathscr G} | \langle f, U(g) \Psi\rangle_{L^2(\mathcal M)} |^2 \d g < \infty 
\end{equation}
for all $f\in L^2(\mathcal M)$. The integral is taken with respect to (any) Haar measure on $\mathscr G$.
Equation \eqref{eq:admCond} implies that the associated wavelet transform
\begin{equation}
  \mathscr W f (g)  = \langle f , U(g)\Psi \rangle_{L^2(\mathcal M)}
\end{equation}
is a bounded operator $\mathscr W: L^2(\mathcal M)\to L^2(\mathscr G)$, while the irreducibility of $U$ ensures invertibility of $\mathscr W$ on its range. To see this, we take an  $f\in\ker \mathscr W$. Then $f\perp\spann\{U(g)\Psi:g\in\mathscr G\}$. Since $U$ is unitary, the orthocomplement of an invariant subspace is invariant and by irreducibility of $U$ this implies $f=0$.

This approach has several disadvantages, which becomes apparent in particular for spheres. Because $L^2(\mathcal M)$ is infinite dimensional, no compact group admits an irreducible unitary representation of this form. However, compact groups seem natural at least in the situation where $M$ itself is a homogeneous space of a compact group. A related approach was incorporated for spheres by Antoine and Vandergheynst. They used quasi-regular representations of the Lorentz group $\SO(n+1,1)$ on the sphere $\S^n \simeq \SO(n+1,1)/\SO(n,1)$. To overcome the the lack of square-integrability they restricted considerations to so-called admissible sections of the Lorentz group $\SO(n+1,1)$. See \cite{Ferreira}, \cite{Ferreira2} for a detailed exposition on that.

An alternative approach to wavelet transforms on spheres was followed by Coifman-Maggioni in \cite{CoifmanMaggioni} or by Freeden \cite{Freeden}. Classical wavelet theory on $\R^n$ is based on the group  generated by translations and dilations. It is evident what translations are rotations on a sphere (seen as homogeneous space of the rotation group), but there is no canonical choice for dilations. The key idea of diffusive wavelets is to generate dilations from a diffusive semigroup, e.g., from time-evolution of solutions to a heat equation on the homogeneous space. 
The use of compact groups is advantageous in view of the availability of powerful tools like Peter-Weyl theorem  and available classifications of irreducible representations. A related concept based on spectral calculus of the Laplacian on closed manifolds was proposed by Geller, \cite{Geller}. 

Discrete wavelet transforms in such a setting were discussed in \cite{CoifmanMaggioni} and \cite{Coifman2}, where heat evolution is combined with a orthogonalisation procedure to model a multi-resolution analysis within $L^2(\S^3)$. Our approach to continuous wavelet transforms given in the present paper follows the discussion of \cite{Freeden} for spheres, but replacing the heat flow by a more general approximate convolution identity.

In the following sections we present a representation theory based approach to understand diffusive wavelet transforms on compact Lie groups and their homogeneous spaces. The wavelet transforms itself are constructed in terms of kernels of an approximate convolution identity in $L^2(\G)$, $\G$ a compact Lie group, or in $L^2(\G/\H)$ with $\H<\G$ a closed subgroup of $\G$. 

The main structure of the paper is as follows. In Section~\ref{sec1} we recall some of the main ingredients from representation theory, in particular Peter-Weyl theorem and the resulting non-commutative Fourier transform on compact groups. For a more detailed introductory account on the analytical aspects and relations to operator theory and partial differential equations of the group Fourier transform we refer to  \cite{RuzhanskyTurunen}.

Section~\ref{sec2} sketches the main ideas of the approach and discusses the construction of diffusive wavelet transforms on arbitrary compact groups. Starting point is the notion of diffusive approximate identities given in Definition~\ref{df1}, for which the heat kernel provides one interesting and simple example. Later all constructions will be explicitly discussed for this particular case.  The concept of diffusive approximate identity is closely related to the scaling function within the classical discrete wavelet approach.

The constructed wavelet transform of $f\in L^2(\G)$ has the usual structure of an inner product between $f$ and a translated wavelet $\psi_\rho$, $\rho$ being a scaling parameter, and also the corresponding inversion formula has the familiar form of a superposition of translations of $\psi_\rho$ integrated over all scales, see formulas \eqref{eq:Wtransf} and \eqref{eq:Wtransf-1}. 

A projection methods transfers these wavelet transforms to homogeneous spaces of $\G$ by identifying them as quotients $\G/\H$ for a suitable subgroup $\H<\G$. On the other hand, the wavelet transform can be intrinsically defined on the homogeneous space itself. Section~\ref{sec3} provides two approaches to wavelet transforms on homogeneous spaces, one giving functions parametrised by points in the space and one parametrised by elements of the group. They extend existing constructions on spheres due to \cite{Freeden} or \cite{BernsteinEbert}.  As one essential ingredient for our approach we include a discussion of several convolution-like products on homogeneous spaces and their relation to the group Fourier transform.

Finally, Section~\ref{sec4} contains concrete examples of transforms arising from our methods. In particular, we discuss the unit circle, the group of unit quaternions and spheres as homogeneous spaces of $\SO(n+1)$. Furthermore, real projective space can be considered analogously to spheres. 

\section{Prerequisites: Representation theory of compact Lie groups}\label{sec1}
We start by recollecting some basic notation. Let $\G$ be a {\em compact} Lie group. A unitary representation of $\G$ is a continuous group homomorphism $\pi:\G\to\U(d_\pi)$ of $\G$ into the group of unitary matrices of a certain dimension $d_\pi$. Such a representation is irreducible if $\pi(g)A=A\pi(g)$ for all $g\in\G$ and some $A\in\C^{d_\pi\times d_\pi}$ implies $A=c\I$ is a multiple of the identity. Equivalently, $\C^{d_\pi}$ does not have non-trivial $\pi$-invariant subspaces $V\subset\C^{d_\pi}$ with $\pi(g)V\subset V$ for all $g\in\G$. We call two representations $\pi_1$ and $\pi_2$ equivalent, if there exists an invertible matrix $A$ with $\pi_1(g)A=A\pi_2(g)$ for all $g\in\G$. 

Let $\widehat\G$ denote the set of all equivalence classes of irreducible representations. Then this set parameterises an orthogonal decomposition of $L^2(\G)$, where $\d\mu_\G$ is the normalised Haar measure on $\G$. 

\begin{thm}[{Peter-Weyl,\cite[Chapter I.4.3]{Vilenkin}}] Let $\G$ be a compact Lie group. Then the following statements are true.
\begin{enumerate}
\item
Denote $\mathcal H_{\pi} =\{g\mapsto \trace(\pi(g)A) : A\in\C^{d_\pi\times d_\pi}\}$. Then the Hilbert space
$L^2(\G)$ decomposes into the orthogonal direct sum
\begin{equation}
   L^2(\G) = \bigoplus_{\pi\in\widehat\G} \mathcal H_\pi.
\end{equation}
\item
For each irreducible representation $\pi\in\widehat\G$ the orthogonal projection $L^2(\G)\to\mathcal H_\pi$ is given by
\begin{equation}
   \phi \mapsto d_\pi \int_\G \phi(h)\chi_\pi (h^{-1}g)\d\mu_\G(h) = d_\pi \phi*\chi_\pi
\end{equation}  
in terms of the character $\chi_\pi(g)=\trace\pi(g)$ of the representation.
\end{enumerate}
\end{thm}

We will denote the matrix $A$ in $\phi*\chi_\pi =\trace(\pi(g)A)$ as Fourier coefficient $\hat\phi(\pi)$ of $\phi$ at the irreducible representation $\pi$. A short calculation yields
\begin{equation}\label{eq:Fdef1}
   \hat\phi(\pi) = \int_\G \phi(g) \pi^*(g)\d\mu_\G(g)
\end{equation}
and therefore
\begin{equation}\label{eq:Fdef2}
   \phi(g) = \sum_{\pi\in\widehat\G} d_\pi \trace(\pi(g) \hat\phi(\pi)).
\end{equation}
If we denote by $\|A\|_{HS}^2 = \trace(A^*A)$ the Frobenius or Hilbert-Schmidt norm of a matrix $A$, then the following
version of Parseval identity holds true.

\begin{cor}[Parseval-identity]
Let $\phi\in L^2(\G)$. Then the matrix-valued Fourier coefficients $\hat\phi(\pi)\in\C^{d_\pi\times d_\pi}$ satisfy
\begin{equation}
   \|\phi\|^2 = \sum_{\pi\in\widehat\G} d_\pi \|\hat\phi(\pi)\|_{HS}^2.
\end{equation}
\end{cor}

The group structure of $\G$ defines an involution of functions on $\G$. We denote it as $\check\phi (g)= \overline{\phi(g^{-1})}$.
On the level of Fourier coefficients this corresponds to taking adjoints of the matrix-valued Fourier coefficients.

On the group $\G$ one defines the convolution of two integrable functions $\phi,\psi\in L^1(\G)$ as
\begin{equation}
  \phi*\psi (g) = \int_\G \phi(h) \psi(h^{-1}g)\d\mu_\G(h).
\end{equation}
Since $\phi*\psi\in L^1(\G)$, the Fourier coefficients are well-defined and they satisfy

\begin{lem}[Convolution theorem on $\G$]
Let $\phi,\psi\in L^1(\G)$ then $\phi*\psi\in L^1(\G)$ and 
\begin{equation}
\widehat{\phi*\psi}(\pi)=\hat\psi(\pi)\hat\phi(\pi).
\end{equation}
\end{lem} 
\begin{proof}
We consider only the convolution formula. As short calculation yields
\begin{align*}
\widehat{\phi*\psi}(\pi) &= \iint \phi(h)\psi(h^{-1}g) \pi^*(g) \d\mu_\G(h)\d\mu_\G(g)\\
&=  \iint \phi(h)\psi(h^{-1}g) \pi^*(h^{-1}g) \pi^*(h) \d\mu_\G(h)\d\mu_\G(g)=\hat\psi(\pi)\hat\phi(\pi).
\end{align*}
\end{proof}

The group structure gives rise to left and right translations, $T_g : \phi \mapsto \phi(g\cdot)$ and $T^g: \phi\mapsto \phi(\cdot g^{-1})$ of functions on the group. For completeness, we remind the reader to the formulae
\begin{equation}
   \widehat{T_g\phi} (\pi) = \hat\phi(\pi)\pi(g),\qquad\text{and}\qquad \widehat{T^g\phi} (\pi) = \pi^*(g)\hat\phi(\pi).
\end{equation}
They are a direct consequence of the definition of the group Fourier transform \eqref{eq:Fdef1}, \eqref{eq:Fdef2}. The assignments $g\mapsto T_g$ and $g\mapsto T^g$ are usually referred to as the left respectively right regular representations of $\G$ in $L^2(\G)$. Both representations are unitary, $T_g^*=T_{g^{-1}}$ and $(T^g)^*=T^{(g^{-1})}$.

We need one more ingredient for later considerations. The Laplace-Beltrami operator $\Delta_\G$ on the group $\G$ is bi-invariant, i.e. it commutes with all $T_g$ and $T^g$. Therefore, all of its eigenspaces are bi-invariant subspaces of $L^2(\G)$. As $\mathcal H_{\pi}$ are minimal bi-invariant subspaces, each of them has to be eigenspace of $\Delta_\G$ and we denote the corresponding eigenvalue as $-\lambda_\pi^2$. Hence, we obtain
\begin{equation}
   \Delta_\G \phi = - \sum_{\pi\in\widehat\G} d_\pi \lambda_\pi^2 \trace(\pi(g)\hat\phi(\pi))
\end{equation}
and the solution to the heat equation
\begin{equation}
  (\partial_t - \Delta_\G)u=0,\qquad u(0,\cdot)=\phi 
\end{equation}
is given as convolution with the heat kernel $p_t(g)$ as $u(t,\cdot)=\phi*p_t$ with
\begin{equation}\label{eq:heatK}
   \widehat p_t(\pi) =  \e^{-t\lambda_\pi^2 } \I,\qquad p_t(g) = \sum_{\pi\in\widehat\G} d_\pi \e^{-t\lambda_\pi^2}\chi_\pi(g). 
\end{equation}
We remark that in particular $\phi*p_t \to\phi$ for all $\phi\in L^p(\G)$, $p\in[1,\infty)$. This fact will be crucial for the proof of inversion formul\ae{} later on.

\section{Diffusive wavelets on a compact Lie group $\G$}\label{sec2}
\subsection{General philosophy}
Let $p_t\in L^1(\G)$ be an approximate convolution identity, i.e., we assume that $\phi*p_t \to \phi$ as $t\to0$ for all $\phi\in L^2(\G)$. We will assign families of convolution kernels $\psi_\rho$ and $\Psi_\rho$ to $p_t$ such that
\begin{equation}\label{eq:adm-cond}
  p_t = \int_t^\infty \check\psi_\rho *\Psi_\rho\, \alpha(\rho)\d\rho.
\end{equation}
We assume that both families belong to $L^1(\G)$. Then we consider the following pair of transformations. Let $\phi\in L^2(\G)$. Then we can assign to $\phi$ a two-parameter function $\mathscr W\phi$ depending on points $g\in\G$ and a scale parameter $\rho\in\R_+$, the {\em wavelet transform}
\begin{equation}\label{eq:Wtransf}
   \phi(g) \mapsto  \mathscr W\phi(\rho,g) = \phi*\check\psi_\rho(g) = \int_\G\phi(h)  \check\psi_\rho(h^{-1}g)\d\mu_\G(h)
   = \langle\phi,T_g^*\psi_\rho \rangle,
\end{equation}
and similarly we can invert this formula by convolution with $\Psi_\rho$ and integrating over the scale parameters
\begin{equation}\label{eq:Wtransf-1}
   \phi = \int_{\to0}^\infty   \mathscr W\phi(\rho,\cdot) * \Psi_\rho \alpha(\rho) \d\rho 
   =\phi * \int_{\to0}^\infty \check\psi_\rho*\Psi_\rho \alpha(\rho)\d\rho .
\end{equation}
Here and later on `${}_{\to0}$' as integral limit means that the integral has to be understood in an improper Riemann sense.

This approach works for arbitrary approximate convolution identities $p_t$. Of interest are in particular those for which the operator $*\partial_t p_t$ is positive. Then the corresponding Fourier coefficients are positive matrices and the choice $\psi_\rho=\Psi_\rho$ seems reasonable.
We will later implement this general philosophy for the particular situation where $p_t$ is the heat kernel and where both families coincide.

\begin{df}\label{df1} Let $\widehat\G_+\subset\widehat\G$ be co-finite. 
A family $t\mapsto p_t$ from $C^1(\R_+;L^1(\G))$ will be called {\em diffusive approximate identity} with respect to $\widehat \G_+$ if
it satisfies
\begin{enumerate}
\item $\|\hat p_t(\pi)\| \le C$ uniform in $\pi\in\widehat\G$ and $t\in\R_+$;
\item $\lim_{t\to0} \hat p_t(\pi) = \I$ for all $\pi\in\widehat\G$;
\item $\lim_{t\to\infty}\hat p_t(\pi) = 0$ for all $\pi\in\widehat \G_+$;
\item $-\partial_t \hat p_t(\pi)$ is a positive matrix for all $t\in\R_+$ and $\pi\in\widehat\G_+$.
\end{enumerate}
\end{df}

\begin{rem}
Conditions (1) and (2) of Definition~\ref{df1} imply by Banach-Steinhaus that $\phi*p_t\to\phi$ as $t\to0$ in $L^2(\G)$.
The norm $\|\cdot\|$ in (1) is the Euclidean matrix norm. 
 Condition (1), (3) and (4) yield $p_t|_{\widehat\G_+} = -\int_t^{\to\infty} \partial_\rho p_\rho\d\rho$ as $L^2(\G)$-valued Riemann integral. 
\end{rem}

\begin{rem}
If $p_t$ is a convolution semi-group, $p_t*p_s=p_{t+s}$, the conditions are equivalent to the negativity of its generator $\mathcal L$ on $(\ker\mathcal L)^\perp = \bigoplus_{\pi\in\widehat\G_+}\mathcal H_\pi$.  In particular the heat kernel $p_t$ as defined in \eqref{eq:heatK} is an diffusive approximate identity with respect to $\widehat\G_+ = \{ \pi\in\widehat\G : \lambda_\pi>0 \} = \{\pi\in\widehat\G :\mathcal H_\pi\ne\C\}$.
\end{rem}

\subsection{Diffusive wavelets} For the following we fix a diffusive approximate identity and denote 
\begin{equation}
L^2_0(\G) =  \bigoplus_{\pi\in\widehat\G_+}\mathcal H_\pi.
\end{equation} 
Wavelet transforms will live on that subspace of square integrable functions. Furthermore, we will require the 
admissibility condition \eqref{eq:adm-cond} for the projection of $p_t$ onto this subspace.
 In the situation of the heat kernel the space $L^2_0(\G)$ contains all square integrable functions with vanishing mean. Furthermore, we denote by $\phi\big|_{\widehat\G_+}$ the orthogonal projection of $\phi\in L^2(\G)$ onto $L^2_0(\G)$ and its obvious extension to distributions on $\G$.

\begin{df} Let $p_t$ be a diffusive approximate identity and $\alpha(\rho)>0$ a given weight function. 
A family $\psi_\rho\in L^2_0(\G)$ is called {\em diffusive wavelet family}, if it satisfies the admissibility condition
\begin{equation}\label{eq:adm-cond-2}
  p_t \big|_{\widehat\G_+} = \int_t^\infty  \check \psi_\rho *\psi_\rho \alpha(\rho)\d\rho.
\end{equation}
\end{df}

Equation \eqref{eq:adm-cond-2} can be solved explicitly. Applying the group Fourier transform to both sides yields
\begin{equation}
  \hat p_t(\pi) =  \int_t^\infty \hat\psi_{\rho}(\pi) \hat\psi_\rho^*(\pi) \alpha(\rho)\d\rho,\qquad \forall\pi\in\widehat\G_+,
\end{equation}
and by differentiating both sides we obtain 
\begin{equation}\label{eq:adm-cond-F}
-\partial_t\hat p_t(\pi) =   \hat\psi_\rho(\pi) \hat\psi_\rho^*(\pi)\alpha(\rho).
\end{equation}
Hence, up to a unitary matrix $\hat\psi_\rho(\pi)$ is the positive square root of the matrix $-\partial_t\hat p_t(\pi)$.

\begin{explH} 
Considering the particular example of heat kernel $p_t$ and resulting heat wavelet families, we can go a step further and write \eqref{eq:adm-cond-F} as
\begin{equation}\label{eq:diffWaveletCoeff}
\hat\psi_\rho(\pi) =\frac1{\sqrt{\alpha(\rho)}} \lambda_\pi \e^{-\rho \lambda_\pi^2 /2} \eta_\pi(\rho)
\end{equation}
for any (fixed) choice of  a family  $\eta_\pi(\rho) \in \U(d_\pi)$. This implies
\begin{equation} 
   \psi_\rho = \frac1{\sqrt{\alpha(\rho)}}\sum_{\pi\in\widehat\G} d_\pi \lambda_\pi \e^{-\rho\lambda_\pi^2/2} \trace(\pi(g)\eta_\pi(\rho)),
\end{equation}
The freedom to choose $\eta_\pi(\rho)$ should not be overrated here; since we want to have $\psi_\rho$ concentrated around the identity element of the group the most sensible choice will be $\eta_\pi(\rho)=\mathrm I$.

The weight function $\alpha(\rho)$ can be used to normalise the family $\psi_\rho$. If we want to have a normalised $L^2$-norm,  Parseval identity yields an expression for the corresponding weight $\alpha(\rho)$
\begin{equation}
  \alpha(\rho) = \sum_{\pi\in\widehat\G} d_\pi^2 \lambda_\pi^2 \e^{-\rho\lambda_\pi^2} =  - \partial_\rho p_\rho(1)  = -\Delta_\G p_\rho(1) \sim \rho^{-1-\frac12 \dim\G}
\end{equation}
based on $\|\eta_\pi\|_{HS}^2 = \trace\I=d_\pi$. Note, that $p_\rho(1)$ is just the heat trace on $\G$. Another notable possibility is to use $\alpha(\rho)=1$ and $\eta_\pi(\rho)=\eta_\pi$ independent of $\rho$ such that $\psi_\rho = p_{\rho/2}*\psi_0$ for $\psi_0=\sum_\pi d_\pi \lambda_\pi\trace(\pi^*(g)\eta_\pi)\in\mathcal D'(\G)$.
\end{explH}

\begin{thm} \label{thm:unit}
Let $\psi_\rho$ be a diffusive wavelet family. Then the associated wavelet transfrom \eqref{eq:Wtransf} $\mathscr W\phi(\rho,\cdot) = \phi*\check\psi_\rho$ is unitary,
\begin{equation}
 \mathscr W : L^2_0(\G)\to L^2(\R_+\times \G, \alpha(\rho)\d\rho \otimes \d\mu_\G).
 \end{equation}
\end{thm}
\begin{proof}
Let $\phi_1$ and $\phi_2$ be two functions with vanishing mean. Then
\begin{align*}
  \langle \mathscr W\phi_1,\mathscr W\phi_2\rangle &=
  \iint \mathscr W\phi_1(\rho,g) \,\overline{\mathscr W\phi_1(\rho,g)}\d\mu_\G(g)\alpha(\rho)\d\rho\\
  &=\int_{\to0}^\infty \iiint \phi_1(h)\check\psi_\rho(h^{-1}g) \overline{\phi_2(h')\check\psi_\rho(h'^{-1}g)}
\d\mu_\G(h) \d\mu_\G(h')\d\mu_\G(g)  \alpha(\rho)\d\rho\\
  &= \lim_{t\to0} \iint    \phi_1(h) \overline{\phi_2(h') }
   \int_{t}^\infty \int \check\psi_\rho(h^{-1}g) \overline{\check\psi_\rho(h'^{-1}g)}\d\mu_\G(g) \alpha(\rho)\d\rho
  \d\mu_\G(h)\d\mu_\G(h')\\
  &=\lim_{t\to0}\iint    \phi_1(h) \overline{\phi_2(h') } \left(p_t(h^{-1}h')-1\right) \d\mu_\G(h')  \d\mu_\G(h)\\ 
  &=\langle\phi_1,\phi_2\rangle
\end{align*}
 based on $\phi*p_t\to\phi$ in $L^2_0(\G)$ as $t\to0$ and the statement follows.
\end{proof}

\begin{rem}
The functions $\psi_\rho$ and $T_g\psi_{\rho'}$ are not orthogonal in general. A simple calculation yields for heat wavelet families the identity
\begin{equation}
   \langle\psi_\rho,T_g\psi_{\rho'}\rangle = \frac1{\sqrt{\alpha(\rho)\alpha(\rho')}} \sum_{\pi\in\widehat\G} d_\pi \lambda_\pi^2 \e^{-\lambda_\pi^2\frac{\rho+\rho'}2}
   \chi_\pi(g) =- \frac1{\sqrt{\alpha(\rho)\alpha(\rho')}}\Delta_\G p_{\frac{\rho+\rho'}2}(g). 
\end{equation}
\end{rem}

\section{Diffusive wavelets on homogeneous spaces $\G/\H$}\label{sec3}
We can use the wavelet transform on $\G$ to define corresponding transforms on homogeneous spaces of $\G$. A particular example we have in mind are spheres $\S^n$ as homogeneous spaces of $\SO(n+1)$. 

Let for the following $\mathcal X$ be a homogeneous space of $\G$, i.e., we assume that there is a transitive group
action $\G \times \mathcal X\to\mathcal X$ denoted as $(g,x)\mapsto g\cdot x$ and satisfying $g_1\cdot g_2\cdot x = (g_1g_2)\cdot x$. If we {\em fix a base-point} $x_0\in\mathcal X$ functions on $\mathcal X$ have a canonic interpretation as functions on $\G$ via
\begin{equation}
   \tilde \phi(g) = \phi(g\cdot x_0). 
\end{equation}
Note, that this identification depends on the chosen base-point in an essential way. Functions obtained in that way are constant on cosets $g\H$ for $\H = \G_{x_0} = \{h\in\G : h\cdot x_0=x_0\}$ the stabiliser of $x_0$ in $\G$.

We can define a $\G$-invariant measure $\d x$ on $\mathcal X$ by setting
\begin{equation}\label{eq:IntDef}
   \int_{\mathcal X} \phi(x) \d x= \int_\G \tilde\phi(g) \d\mu_\G(g)
\end{equation}
and ask for a wavelet transform within $L^2(\mathcal X)$. 
\begin{lem}
Equation \eqref{eq:IntDef} is independent of the choice of the base-point $x_0$ and defines a $\G$-invariant Radon measure on the homogeneous space $\mathcal X$.
\end{lem}
There are two possible ways to define the wavelet transform. The first one is to apply the previously defined transform to the lifted function on the group $\G$, the second one is to define a genuine transform associated directly to the homogeneous space. We will discuss both approaches in the next subsections. For simplicity we will first treat the 
simpler situation where wavelets on the group are expanded in terms of characters.

\subsection{The naive way}\label{sec:zonal1}
We now apply the wavelet transform to the lifted function $\tilde\phi$ for some $\phi\in L^2(\mathcal X)$. This defines a function on $\R_+\times\G$ via
\begin{equation}
 \mathscr W\tilde \phi(\rho,g) =  \int_\G \tilde \phi(h) \check\psi_\rho(h^{-1}g) \d\mu_\G(h)        =  \int_\G   \phi(h\cdot x_0) \check\psi_\rho(h^{-1}g) \d\mu_\G(h). 
\end{equation}
This definition has one obvious drawback; we would prefer to have a transform living on $\R_+\times\mathcal X$ instead of $\R_+\times\G$. 

\begin{explH}
Heat wavelet families can be defined in terms of characters, i.e. 
\begin{equation}
   \psi_\rho(g) = \frac1{\sqrt{\alpha(\rho)}} \sum_{\pi\in\widehat\G} d_\pi \lambda_\pi \e^{-\lambda_\pi^2\rho/2} \eta_\pi(\rho)
   \chi_\pi(g)
\end{equation}
with a (fixed) choice of a family $\eta_\pi(\rho)\in\C$, $|\eta_\pi(\rho)|=1$. Then $\psi_\rho$ is a central function on $\G$, i.e. it is constant on conjugacy classes. Then the naively defined wavelet transform lives on the homogeneous space $\mathcal X$, 
\begin{align}\label{eq:ZonalWtransf}
  \mathscr W_{\mathcal X} \phi(\rho,x)  &=  \mathscr W \tilde \phi(\rho,g) =   \int_\G   \phi(h\cdot x_0) \check\psi_\rho(h^{-1}g) \d\mu_\G(h)\notag\\& = \int_\G   \phi(h\cdot x) \overline{\psi_\rho(h)} \d\mu_\G(h)
\end{align}
where $g$ is such that $x=g\cdot x_0$.
\end{explH}

The same reasoning works as long as the wavelet family on $\G$ consists of central functions. 
We will refer to this transform as (naive) {\em zonal wavelet transform} on $\mathcal X$. We can write \eqref{eq:ZonalWtransf}
as integral over $\mathcal X$. Indeed, averaging formula \eqref{eq:ZonalWtransf} over $\H$ yields 
\begin{align}
  \mathscr W_{\mathcal X} \phi(\rho,x) & = \int_\G   \phi(h\cdot x_0) \int_\H \overline{\psi_\rho(g^{-1}hh')} \d\mu_\H(h') \d\mu_\G(h)\notag\\
  &=   \int_\G   \phi(h\cdot x_0)  \overline{\mathbb P_{\mathcal X} \psi_\rho(g^{-1}h \cdot x_0) } \d\mu_\G(h) \notag\\
  &= \int_{\mathcal X} \phi(y) 
 \overline{\mathbb P_{\mathcal X}  \psi_\rho(g^{-1}\cdot y)}\d y 
\end{align}
with
\begin{equation}
   \mathbb P_{\mathcal X} \psi_\rho (h\cdot x_0) = \int_\H \psi_\rho(hh') \d\mu_\H(h')
\end{equation}
the $\H$-average of $\psi_\rho$ (which is constant on cosets $h\H$) seen as function on the homogeneous space $\mathcal X$ with fixed base point $x_0$. The measure $\d\mu_\H$ is the (normalised) Haar measure on $\H$.

\begin{cor} Let $\psi_\rho$ be a diffusive wavelet family and assume that $\psi_\rho$ are central functions on $\G$.
Then the associated zonal wavelet transform $\mathscr W_{\mathcal X}$ is unitary $L^2_0(\mathcal X)\to L^2(\R_+\times \mathcal X,\alpha(\rho)\d\rho\otimes\d x)$. 
\end{cor}
\begin{proof}
The calculation is reduced to the one on $\G$. Indeed,
\begin{align*}
   \langle  \mathscr W_{\mathcal X}\phi_1, \mathscr W_{\mathcal X}\phi_2\rangle
   &= \int_{\to0}^\infty \int_{\mathcal X}  \mathscr W_{\mathcal X}\phi_1(\rho,x)\overline{\mathscr W_{\mathcal X}\phi_2(\rho,x)}\d x\alpha(\rho)\d\rho\\
   &=\langle  \mathscr W\tilde \phi_1, \mathscr W\tilde \phi_2\rangle
   =\langle\tilde\phi_1,\tilde\phi_2\rangle = \langle\phi_1,\phi_2\rangle
\end{align*}
follows from Theorem~\ref{thm:unit} and \eqref{eq:IntDef}.
\end{proof}

\subsection{Zonal wavelet transforms}\label{sec:zonal}
We want to have a closer look at the structure of our formulas. For this we recall some notation. 
\begin{df}
A function $\phi$ on $\mathcal X$ is {\em zonal} with respect to the base point $x_0$ if it is invariant under the action of the stabiliser $\H$ of $x_0$.
\end{df}
This is the case for $\mathbb P_{\mathcal X}\psi_\rho$ and true in general for any projection of a central function $\chi$ on $\G$
\begin{align}
   \mathbb P_{\mathcal X} \chi(h\cdot x) &=     \mathbb P_{\mathcal X} \chi(hg\cdot x_0) =
 \int_\H \chi(hgh') \d\mu_\H(h') \notag\\&=    \int_\H \chi(gh'h) \d\mu_\H(h') =    \mathbb P_{\mathcal X} 	\chi(x).
\end{align}

To avoid a too cumbersome notation, we will simply and write $\hat\phi$ instead of $\hat{\tilde\phi}$ for Fourier coefficients of functions on $\mathcal X$ (meaning that we identify $\mathcal X$ for fixed base-point $x_0$ with $\G/\H$).  

\begin{prop}
There exists a family of orthogonal projections $\hat\H(\pi)\in\C^{d_\pi\times d_\pi}$ depending on the stabiliser $\H$ such that 
\begin{equation}
 \widehat{\mathbb P_{\mathcal X} \phi}(\pi)
 = \hat\H(\pi) \hat\phi(\pi).
\end{equation}
\end{prop}
\begin{proof}
A short calculation yields
\begin{equation}
\hat\H(\pi) =  \int_\H  \pi(h) \d\mu_\H(h),
\end{equation}
such that $\hat\H(\pi)\hat\H(\pi)=\hat\H(\pi)$ and $\hat\H^*(\pi)=\hat\H(\pi)$.
\end{proof}

\begin{cor}
A function $\phi\in L^1(\mathcal X)$ is zonal if and only if $\hat\H(\pi)\hat\phi(\pi)=\hat\phi(\pi)\hat\H(\pi)$ for all $\pi\in\widehat\G$. This is true, if and only if  there exists a function $\check\phi\in L^2(\mathcal X)$ with $\check{\tilde\phi}=\tilde{\check\phi}$.
\end{cor}

Since $\hat\H(\pi)$ is an orthogonal projector in $\C^{d_\pi}$ we can choose a basis of $\C^{d_\pi}$ in such a way that $\hat\H(\pi)=\bdiag(\I_{r_\pi},0)$ is block-diagonal, $r_\pi=\rank\hat\H(\pi)$. Then, for $\phi\in L^1(\mathcal X)$ 
Fourier coefficients satisfy $(\mathrm I-\hat\H(\pi))\hat\phi(\pi)=0$, which means, it is a matrix having non-zero entries in the first $r_\pi$ rows only. Similarly, $\phi\in L^1(\mathcal X)$ is zonal if its Fourier coefficients have entries only in the upper-left $r_\pi\times r_\pi$ block. 

Convolutions of functions defined on $\mathcal X$ should be related to multiplications of Fourier coefficients. The underlying block structure of Fourier coefficients makes it advantageous to combine multiplications with taking adjoints. This leads to the introduction of several convolution like products for functions; we collect them in the following definition. 
\begin{df}  Let $\phi,\psi\in L^1(\mathcal X)$. Then we define
\begin{enumerate}
\item the {\em group convolution}
\begin{equation}
    \phi * \psi (x) = \int_\G \phi(g\cdot x_0) \psi(g^{-1}\cdot x) \d\mu_{\G}(g)  \in L^1(\mathcal X);
\end{equation}
\item the {\em $\bullet$-product}
\begin{equation}
    \phi \bullet \psi (g) = \int_{\mathcal X} \phi(x) \overline{\psi(g^{-1}\cdot x)} \d x = \langle  \phi,T_g^* \psi\rangle \in L^1(\G);
\end{equation}
\item the {\em zonal product} 
\begin{equation}
    \phi \,{\hat\bullet}\, \psi (x) = \int_{\G} \overline{\phi(g\cdot x_0)}{\psi(g\cdot x)} \d\mu_\G(g) \in L^1(\mathcal X).
\end{equation}
\end{enumerate}\end{df}

In order to justify the above definitions and to relate them to Fourier coefficients we collect some properties of these convolution like products in the following proposition.

\begin{prop} Let $\phi,\psi\in L^1(\mathcal X)$.
\begin{enumerate}
\item $\widetilde{\phi*\psi} = \tilde\phi*\tilde\psi$ and thus $\widehat{\phi*\psi}(\pi)=\hat\psi(\pi) \hat\phi(\pi)$. 
\item $\phi\bullet\psi = \tilde\phi * \check{\tilde\psi} $ and thus $\widehat{{\phi\bullet\psi}}(\pi) = \widehat{\psi}^*(\pi)\widehat{\phi}(\pi)$.
\item If $\psi$ is zonal with respect to $x_0$ then $\phi\bullet\psi$ is constant on cosets $g\H$ and thus defines a function on $\mathcal X$.
\item $\phi\,\hat\bullet\,\psi$ is zonal with respect to $x_0$.
\item ${\widetilde{\phi\,\hat\bullet\,\psi}} =\check{\tilde\phi}*{\tilde\psi}$ and
 $\widehat{{\phi\,\hat\bullet\,\psi}}(\pi) = \widehat{\psi}(\pi)\widehat{\phi}^*(\pi)$.
\end{enumerate}
\end{prop}
\begin{proof}
(1) evident. (2) It suffices to consider the first formula which follows directly from \eqref{eq:IntDef}. (3) We calculate
$$ \phi\bullet\psi(gh) = \int_{\mathcal X} \phi(x)\overline{\psi(h^{-1}g^{-1}\cdot x)} \d x= \int_{\mathcal X} \phi(x)\overline{\psi(g^{-1}\cdot x)} \d x= \phi\bullet\psi(g) $$
based on the zonality of $\psi$. (4) For $h\in\H=\G_{x_0}$ we obtain
$$ \phi\,{\hat\bullet}\,\psi(h\cdot x) = \int_\G \overline{\phi(g\cdot x_0)}\psi(gh\cdot x)\d\mu_\G(g)=
 \int_\G \overline{\phi(gh^{-1}\cdot x_0)}\psi(g\cdot x)\d\mu_\G(g) = \phi\,{\hat\bullet}\,\psi(x) $$
by the right invariance of $\mu_\G$. (5) It suffices to show the first formula. Substituting $g^{-1}$ for $g$ and using  \eqref{eq:IntDef} implies the assertion.
\end{proof}

\begin{cor}\label{cor:zonalConv} Let $\phi,\psi\in L^1(\mathcal X)$. 
  If $\psi$ is zonal then $\phi*\psi=\phi\bullet\check\psi$.
  Similarly, if $\phi$ is zonal then $\phi*\psi = \check\phi\,\hat\bullet\,\psi$.
\end{cor}

\begin{cor}
The (naive) zonal wavelet transform is given by a $\bullet$-product
\begin{equation}
   \mathscr W_{\mathcal X}\phi(\rho,x) = \phi \bullet \mathbb P_{\mathcal X}\psi_\rho
\end{equation}
with inversion  
\begin{equation}
   \phi(x) =    \int_{\to0}^\infty \mathscr W_{\mathcal X}\phi(\rho,\cdot) \bullet  \mathbb P_{\mathcal X}\check\psi_\rho  \,\alpha(\rho)\d\rho.
\end{equation}
\end{cor}

The $\bullet$-product of zonal functions is not associative. However, a short computation (e.g., on the Fourier side) yields $(\phi\bullet\psi)\bullet\chi = \phi\bullet(\chi*\psi)$. Therefore, the admissibility condition implies or rather reads as
\begin{equation}
\delta_{x_0}  \big|_{\widehat\G_+}= \int_{\to0}^\infty \mathbb P_{\mathcal X} \check\psi_\rho * \mathbb P_{\mathcal X}\psi_\rho  \,\alpha(\rho)\d\rho.
\end{equation}
In fact, this can be used as starting point to define {\em zonal wavelet transforms} on $\mathcal X$.

\begin{df}\label{df3}
Let $p_t$ be a diffusive approximate identity and let $\alpha(\rho)\ge0$ be a given weight function.
A family $\psi_\rho\in L^2(\mathcal X)$ is called {\em zonal diffusive wavelet family} if
\begin{enumerate}
\item $\psi_\rho$ is zonal with respect to $x_0$,
\item the admissibility condition
\begin{equation}\label{eq:adm-cond-3}
   p_t^{\mathcal X}(x) \big|_{\widehat\G_+} = \int_{t}^\infty \check\psi_\rho *\psi_\rho(x)\; \alpha(\rho)\d\rho
\end{equation}
is satisfied, where 
\begin{equation}
 p_t^{\mathcal X}(g\cdot x_0) =  \int_\H\int_\H p_t( h_1 g h_2) \d\mu_\H(h_1)\d\mu_\H(h_2) 
\end{equation} 
 is the zonal average of $p_t$ on $\mathcal X$.
\end{enumerate}
\end{df}

We associate to this family the wavelet transform $\mathscr W_{\mathcal X} \phi(\rho,x) = \phi\bullet \psi_\rho(x)$ with inversion given as
\begin{equation}
    \phi = \int_{\to0}^\infty \mathscr W_{\mathcal X}\phi(\rho,\cdot)\bullet\check\psi_\rho\,\alpha(\rho)\d\rho
\end{equation}
for all functions $\phi\in L^2_0(\mathcal X)$.

\begin{thm}\label{thm:unit2}
The zonal diffusive wavelet transform is unitary $L^2_0(\mathcal X)\to L^2(\R_+\times \mathcal X, \alpha(\rho)\d\rho\otimes\d x)$.
\end{thm}
\begin{proof}
This can be shown along similar lines to the proof of Theorem~\ref{thm:unit}. Alternatively, one can use Plancherel's theorem and show that the corresponding transform in Fourier space is unitary. We follow that line of thought and write
\begin{align*}
  \langle\phi_1,\phi_2\rangle &= \sum_{\pi\in\widehat\G_+} d_\pi \trace(\hat\phi_2^*(\pi)\hat\phi_1(\pi))= \sum_{\pi\in\widehat\G_+} d_\pi \trace(\hat\phi_2^*(\pi)\hat\H(\pi)\hat\phi_1(\pi))\\
  &=\int_0^\infty \sum_{\pi\in\widehat\G_+} d_\pi \trace(\hat\phi_2^*(\pi) \hat\psi_\rho(\pi) \hat\psi_\rho^*(\pi)\hat\phi_1(\pi))  \alpha(\rho)\d\rho\\
  & =\langle\mathscr W_{\mathcal X}\phi_1,\mathscr W_{\mathcal X}\phi_2\rangle
\end{align*}
based on  \eqref{eq:adm-cond-3} in the form of $\int_0^\infty \hat\psi_\rho(\pi) \hat\psi_\rho^*(\pi)\alpha(\rho)\d\rho = \hat\H(\pi)$ for $\pi\in\widehat\G_+$.
\end{proof}

\subsection{Nonzonal wavelets}\label{sec:nonzonal}
Now we will consider the situation where $\psi_\rho\in L^2(\mathcal X)$ is an arbitrary (non-zonal) family. Then
\begin{equation}
   \mathscr W\phi (\rho,g) = \phi \bullet \psi_\rho (g) = \langle \phi,T_g^*\psi_\rho\rangle
\end{equation}
lives on $\G$ rather than on $\mathcal X$. We aim for an inversion formula of the kind
\begin{equation}
   \phi(x) = \mathbb P_{\mathcal X} \int_{\to0}^\infty  \mathscr W\phi (\rho,\cdot) * \tilde \Psi_\rho (g) \,\alpha(\rho)\d\rho
\end{equation}
with a second family $\Psi_\rho\in L^2(\mathcal X)$. By a short computation we see that $(\phi\bullet\psi)*\tilde \chi = \phi\bullet(\chi\,\hat\bullet\,\psi)$ for $\phi,\psi,\chi\in L^1(\mathcal X)$ and this motivates us to define

\begin{df}\label{df4} 
Let $p_t$ be a diffusive approximate identity and $\alpha(\rho)\ge0$ be a given weight function.
A family $\psi_\rho\in L^2(\mathcal X)$ is called (non-zonal) {\em diffusive wavelet family} if the admissibility condition
\begin{equation}\label{adm-cond-3}
   p_t ^{\mathcal X}(x)\big|_{\widehat\G_+} = \int_t^\infty \psi_\rho\,\hat\bullet\, \psi_\rho(x)\, \alpha(\rho)\d\rho
\end{equation}
is satisfied.
\end{df}
We associate to this family the wavelet transform $\mathscr W_{\mathcal X} \phi(\rho,g) = \phi\bullet \psi_\rho(g)$ with inverse given as
\begin{equation}
    \tilde\phi = \int_{\to0}^\infty \mathscr W_{\mathcal X}\phi(\rho,\cdot)*\tilde \psi_\rho\,\alpha(\rho)\d\rho
\end{equation}
for all $\phi\in L^2_0(\mathcal X)$.

\begin{rem}
If $\psi_\rho$ is zonal, $\check\psi_\rho*\psi_\rho = \psi_\rho\,\hat\bullet\,\psi_\rho$ and $ \mathscr W_{\mathcal X}\phi(\rho,\cdot)*\tilde \psi_\rho =  \mathscr W_{\mathcal X}\phi(\rho,\cdot)\bullet\check \psi_\rho$ by Corollary~\ref{cor:zonalConv}. Hence Definition~\ref{df4} extends Definition~\ref{df3}.
\end{rem}

The proof given to Theorem~\ref{thm:unit2} transfers to the non-zonal situation. Hence, we obtain

\begin{cor}
The (non-zonal) diffusive wavelet transform is unitary $L^2_0(\mathcal X)\to L^2(\R_+\times \G, \alpha(\rho)\d\rho\otimes\d\mu_\G)$.
\end{cor}

\begin{rem}
We want to emphasise the compatibility of the wavelet transforms on $\mathcal X$ and on $\mathscr G$. Therefore,
we look at the transformations in Fourier space. The admissibility equation reads in both cases as an integral of
$\hat\psi_\rho(\pi) \hat\psi_\rho^*(\pi)$. If $\psi_\rho$ is a wavelet on $\mathcal X$ the result are Fourier coefficients
of a zonal function.

The wavelet transform on $\mathscr  G$ as well as on $\mathcal X$ can be written in Fourier space as
\begin{equation}
\widehat{ \mathscr W \phi}(\rho,\pi)=\hat \psi^*_\rho(\pi) \hat \phi(\pi).
\end{equation}
On $\mathcal X$ this gives in general a Fourier coefficient for a function on
$\mathscr G$ (i.e., a full matrix). $\widehat{\mathscr W \phi}$ is the Fourier coefficient of a function on $\mathcal
X$ if and only if $\hat\psi_\rho$ is Fourier coefficient of a zonal function. 
In both cases the reconstruction is given as integral 
\begin{equation}
\hat \phi (\pi)=\int_{\to0}^\infty\hat\psi_\rho(\pi) \widehat{\mathscr W \phi}(\rho,\pi)\alpha(\rho)\d\rho
=  \int^\infty_{\to0} \hat \psi_\rho(\pi) \hat\psi_\rho^*(\pi) \alpha(\rho)\d\rho \hat f(\pi).
\end{equation}
\end{rem}

\begin{explH}
We can describe all admissible non-zonal heat wavelet families by looking at their Fourier coefficients. Formula~\eqref{adm-cond-3} implies $\hat\psi_\rho(\pi)=0$ for $\lambda_\pi=0$ and 
\begin{equation}\label{adm-cond-3'}
 \lambda_\pi^2 \hat\H(\pi) \e^{-\lambda_\pi^2 \rho} =  \hat \psi_\rho(\pi) \hat\psi_\rho^*(\pi) \alpha(\rho),\qquad \lambda_\pi>0.
\end{equation}
We assume that we have chosen the basis  in $\C^{d_\pi}$ in such a way that $\hat\H(\pi)=\bdiag(\I_{r_\pi},0)$, $r_\pi=\rank\hat\H(\pi)$. Then, $\hat\psi_\rho(\pi)$ has non-zero entries in the first $r_\pi$ rows only and condition \eqref{adm-cond-3'} implies for these entries that (i) the rows are pair-wise orthogonal 
\begin{equation}\label{eq:adm-cond-3mass1}
   \sum_{k=1}^{d_\pi}(\hat \psi_\rho(\pi))_{ik} \overline{(\hat \psi_\rho(\pi))_{jk}} = 0,\qquad i\ne j,\quad i,j=1,\ldots,r_\pi
\end{equation}
and (ii) of prescribed length
\begin{equation}\label{eq:adm-cond-3mass2}
   \sum_{k=1}^{d_\pi}|(\hat \psi_\rho(\pi))_{ik}|^2  = \frac1{\alpha(\rho)} \lambda_\pi^2 \e^{-\lambda_\pi^2 \rho},\qquad i=1,\ldots,r_\pi .
\end{equation}
Of particular interest are situations where $r_\pi\in\{0,1\}$. Then \eqref{eq:adm-cond-3mass1} disappears and 
\eqref{eq:adm-cond-3mass2} reduces to one condition. This corresponds to the assumptions used in \cite{BernsteinEbert} and \cite{EBCK}.
\end{explH}

We want to remind the reader of the following definition, which is particularly helpful when considering concrete geometric situation as we do in Sections~\ref{sec4.2} and~\ref{sec4.3}. 
\begin{df}
The subgroup $\H$ is called massive, if $\rank\hat\H(\pi)\le1$ for all $\pi\in\widehat\G$. Furthermore, an irreducible representation $\pi\in\widehat\G$ is called class-1 representation with respect to the subgroup $\H$, if $\rank\hat\H(\pi)\ge1$.
\end{df}

\section{Examples}\label{sec4}

\subsection{Tori and $\vartheta$-functions} First, we consider the simplest possible situation. Let $\G=\mathbb T=\U(1)\subset\C$ be the set of unimodular complex numbers. Then $\widehat{\mathbb T}=\Z$ and to each $k\in\Z$ corresponds an irreducible representation $z\mapsto z^k$. The corresponding Fourier series are just Laurent series or the usual Fourier series if we write $z=\exp(2\pi\i\theta)$. The Laplacian on $\mathbb T$ corresponds in this notation to $\partial_\theta^2$ and has eigenvalues $-4\pi^2k^2$.

In this notation, the heat kernel on $\mathbb T$ is given by
\begin{equation}
   p_t(z) = \sum_{k\in\Z} \e^{-4\pi^2k^2t }z^k =1 +2 \sum_{k=1}^\infty  \e^{-4\pi^2k^2t }\cos(2k\pi \theta) = \vartheta_3(\pi\theta, \e^{-4\pi^2t})
\end{equation}
in terms of Jacobi's $\vartheta_3$-function, cf. \cite[Chapter XXI]{Wittaker:1922}. Similarly, one obtains
\begin{equation}
 \sqrt{\alpha(\rho)} \psi_\rho(z) = -2\pi\i  \sum_{k\in\Z}  k e^{-2\pi^2k^2\rho}z^k = 4\pi \sum_{k=1}^\infty  k \e^{-2\pi^2k^2\rho}\sin(2k\pi \theta) =  -\partial_\theta  \vartheta_3(\pi\theta,\e^{-2\pi^2\rho})
\end{equation}
(choosing $\eta_k=-\i\,\mathrm{sign}\, k$ to simplify the notation). The corresponding wavelet transform of a function 
$f\in L^2[0,1]\simeq L^2(\mathbb T)$ with normalisation $\alpha(\rho)=1$ is
\begin{align}
   \mathscr W \phi(\rho,\theta) &= \int_0^1 f(\tau) \partial_\tau \vartheta_3\left(\pi (\tau-\theta), \e^{-2\pi^2\rho}\right)\d\tau \notag\\
   &=   \int_0^1 f'(\theta-\tau) \vartheta_3(\pi\tau ,\e^{-2\pi^2\rho})\d\tau
\end{align}
with inversion formula 
\begin{equation}
   \phi(\theta) =\int \phi(\tau)\d\tau - \int_{\to0}^\infty \int_0^1 \mathscr W\phi(\rho,\theta-\tau) \, \partial_\tau\vartheta_3(\pi\tau,\e^{-2\pi^2\rho} )\d\tau\d\rho.
\end{equation}
The wavelet transform $\mathscr W\phi(\rho,\theta)$ describes for {\em small} $\rho$ the `high-frequency part' of $\phi$ localised near the point $\theta$. 

\begin{figure}
\includegraphics[width=11cm]{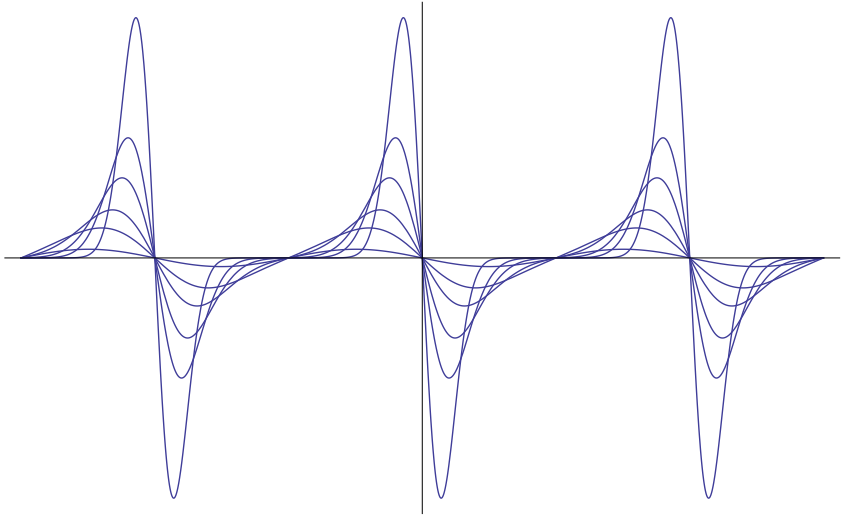}
\caption{The toroidal family $\vartheta_3'(\pi\theta,\e^{-2\pi^2\rho})$ for $-\frac32\le\theta\le\frac32$ and scale parameters $\rho\in\{ 0.005, 0.01, 0.015, 0.025, 0.04, 0.1\}$.}
\label{fig1}
\end{figure}

This $\vartheta$-transformation generalises to higher dimensional tori in the obvious way. We do not give further formul\ae{}, but conclude this example with some pictures of the family $\psi_\rho$ for different $\rho$ depicted in Figure~\ref{fig1}.

\subsection{The group $\S^3$}\label{sec4.2} The sphere $\S^3$ can be viewed as group of unit quaternions, $\S^3\subset\mathbb H$, or alternatively, as group of unitary $2\times2$ matrices $\mathrm{SU}\,(2)$ (see, e.g., \cite{RuzhanskyTurunen2} for that approach). We decide for quaternions $x = x_0+ \mathrm i x_1 + \mathrm j x_2 + \mathrm k x_3\in\mathbb H$  with $x_i\in\R$, $\mathrm i^2=\mathrm j^2=\mathrm k^2=\mathrm{ijk}=-1$. The unitary dual $\widehat{\S^3}$ can be identified with $\mathbb N_0$; for each $n\in\N_0$ there is a representation
\begin{equation}
  t_n : \S^3 \to \U(\Pi_n), \qquad t_n(x) p_n = p_n(x^{-1}\cdot),
\end{equation}
on the space $\Pi_n\subset \C[x_0,\ldots, x_3]$ of all polynomials in four variables which are homogenenous
of degree $n$ and $L^2(\S^3)$ perpendicular to $\Pi_{n'}$, $n'<n$. As $\bigcup_n\Pi_n$ is dense in $L^2(\S^3)$ we obtain all representations in this way. Furthermore, $\Pi_n$ splits into $d_n$ minimal invariant subspaces of dimension $d_n$ corresponding to the splitting of $t_n$ into irreducible components.

\begin{figure}
\includegraphics[width=7cm]{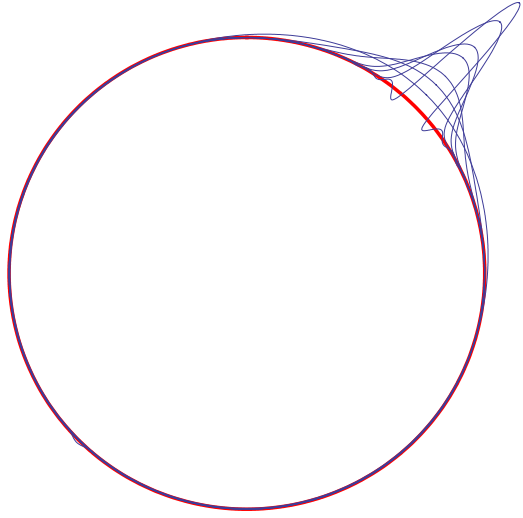}
\caption{The constructed heat wavelet family on $\S^3$. For the figure we have chosen $\alpha(\rho)\sim\rho^{-1}$ and $\rho\in\{0.001,  0.003, 0.005, 0.007, 0.01,   0.02,  0.03\}$.}
\label{fig2}
\end{figure}

An orthonormal basis of $\Pi_n$ is given by spherical harmonics, which occur as matrix coefficients of our representation $\{\mathcal T_n^{ij}(x): i,j=1,\ldots d_n\}$, $d_n=n+1$. As the maximal torus of $\S^3$ is of dimension 1 (and we may fix it as $\S^3\cap\C=\{x\in\S^3:x_2=x_3=0\}$) class functions are functions of just one variable $x_0=\mathrm{Sc}\,x$. This is seen easily from
$y^{-1} x y = \mathrm{Sc}\, x + y^{-1}(\mathrm{Vec}\,x)y$.  Furthermore, such functions are essentially unique within $\Pi_n$; they are multiples of the Gegenbauer polynomial $\mathcal C_n^1 (x_0)$. Hence, using $\mathcal C^1_n(1)=n+1=d_n$ we obtain
\begin{equation}
   \trace t_n (x) = \sum_{j=1}^{d_n} \mathcal T_n^{jj}(x) =  \mathcal C_n^1 (x_0),\qquad x_0=\mathrm{Sc}\, x.
\end{equation}
We normalise the measure on $\S^3$ to be the Haar measure (which changes the Laplacian from the usual one). Then
$\lambda_n^2= (2\pi^2)^{2/3} n(n+2)$ and we obtain the corresponding heat kernel on $\S^3$ as
\begin{equation}
   p_t(x) = \sum_{n=0}^\infty (n+1) \e^{-\lambda_n^2t} \mathcal C_n^1 (x_0),\qquad x_0=\mathrm{Sc}\, x,
\end{equation}
such that the corresponding heat wavelet family is given as
\begin{equation}\label{eq:S3-wavelet}
   \psi_\rho(x) = \frac1{\sqrt{\alpha(\rho)}}\sum_{n=1}^\infty(n+1) {\lambda_n}  \e^{- \lambda_n^2\rho/2} \mathcal C_n^1 (x_0),\qquad x_0=\mathrm{Sc}\, x.
\end{equation}
Figure~\ref{fig2} depicts $\psi_\rho(x)$ on $\S^3\cap \C$ for different $\rho$.

Some spherical 3-manifolds can be treated similarly to $\S^3$ in the previous section. Let for this $\Gamma$ be a {\em finite} subgroup of $\S^3\subseteq \mathbb H$ and $\mathcal X = \S^3 / \Gamma$.  
Interpreting all integrals over $\Gamma$ as averages, we can retrace the considerations of Sections~\ref{sec:zonal1}
to \ref{sec:nonzonal}. A zonal heat wavelet family on $\S^3/\Gamma$ is therefore given by
\begin{equation}
  \psi_\rho(x) = \frac1{\sqrt{\alpha(\rho)}}\sum_{n=1}^\infty (n+1) {\lambda_n}  \e^{- \lambda_n^2\rho/2} 
\left(  \frac1{|\Gamma|} \sum_{\gamma\in\Gamma} \mathcal C_n^1 (\mathrm{Sc}\, x\gamma)\right).
\end{equation}
Main ingredient to understand zonal and non-zonal transforms on the Fourier level is the family of matrices
\begin{equation}
   (\widehat\Gamma(t_n))_{ij} = \frac1{|\Gamma|} \sum_{\gamma\in\Gamma}   \mathcal T_n^{ij}(\gamma).
\end{equation}
As example we consider the lens spaces $L(p,1)$. Then $\Gamma=\{\omega\in\S^3\cap\C : \omega^p=1\}\simeq \Z/p\Z$ and a simple calculation (based on \cite[10.8(17)]{Erdely}) yields
\begin{align}
  \rank \widehat\Gamma(t_n) &= \trace \widehat\Gamma(t_n) = \frac1p \sum_{k=1}^p \mathcal C^1_n\big(\cos \frac{2 k \pi}p\big)
  = \frac1p\sum_{k=1}^p \sum_{m=0}^n  \cos \frac{2(n-2m)k \pi}{p}  \notag\\
  &=\# \{ m\in\{0,\ldots,n\} : n-2m \equiv 0 ( p ) \}.
\end{align}
In particular for the projective space $\mathbb{RP}^3=L(2,1)$ it follows that $\widehat\Gamma(t_n) =\mathrm I$ for $n$ even and $\widehat\Gamma(t_n)=0$ for $n$ odd. The zonal heat wavelet transform on projective space is thus obtained by ignoring all odd Fourier coefficients in \eqref{eq:S3-wavelet}.

\subsection{$\SO(n+1)$ and $\S^n$}\label{sec4.3}
Part of the computations for this example are taken from \cite{EBCK}. We will just recall some notation, for more details see \cite{Gradshteyn} or \cite{Vilenkin}.  
As usual we denote by $\mathcal Y_{k}^i\in C^{\infty}(\S^n)$, $k\in\N_0$, $i=1,\ldots, d_k(n)$ the orthonormal system of spherical harmonics on $\S^n$ (normalised with respect to Lebesgue measure on $\S^n$). Their span is invariant under 
the action of $\SO(n+1)$ and therefore it makes sense to define
the Wigner polynomials $\mathcal T_k^{ij}(A)$, $A\in\SO(n+1)$,
\begin{equation}\label{eq:wigner}
 \mathcal Y_k^i(A^{-1}\xi) = \sum_{j=1}^{d_k(n)} \mathcal T_k^{ij}(A) \mathcal Y_k^j(\xi),\qquad
 \mathcal T_k^{ij}(A) = \int_{\S^n}  \mathcal Y_k^i(A^{-1}\xi) \overline{\mathcal Y_k^j(\xi)}\d\xi
\end{equation}
as matrix coefficients of the corresponding representation. The matrix-valued function $A\mapsto \mathcal T_k(A)$ forms an irreducible unitary class-1 representation of $\SO(n+1)$ on $\S^n\simeq\SO(n+1)/\SO(n)$. 

\begin{rem}
The Wigner polynomials do not give all irreducible unitary representations of $\SO(n+1)$. As we are mainly interested in $\S^n$ here, we do not need to care about that. Lifts of functions on $\S^n$ to $\SO(n+1)$ have vanishing Fourier coefficients on all other representations.
\end{rem}

\begin{lem}[{\cite[Chapter IX.2]{Vilenkin}}]
$\SO(n)$ is a massive subgroup of $\SO(n+1)$. Furthermore, the family $\mathcal T_k$, $k\in\N_0$, gives up to equivalence all class-1 representations of $\SO(n+1)$ with respect to $\SO(n)$. 
\end{lem}

For the following we fix the `north pole' $\xi_0$ of $\S^n$. 
Then it follows that the set of zonal spherical harmonics (w.r.t. $\xi_0$) is one-dimensional and spanned by the Gegenbauer polynomial $\mathcal C_k^{(n-1)/2}(\xi_0^\top\xi)$. 
This implies several interesting formul\ae{}, in particular 
\begin{lem}[Addition theorem]
For all $\xi,\eta\in\S^n$ and $k\in \N_0$
\begin{equation}\label{eq:4.9}
 \frac{\mathcal C_k^{(n-1)/2} (\xi^\top\eta)}{\mathcal C_k^{(n-1)/2}(1)}
  = \frac{|\S^n|}{d_k(n)} \sum_{i=1}^{d_k(n)} \mathcal Y^i_k(\xi)\overline{\mathcal Y^i_k(\eta)}.
\end{equation}
\end{lem}
\begin{proof}[Sketch of proof]
It suffices to check that the right hand side is zonal with respect to $\xi$, which follows from \eqref{eq:wigner} with $A$ in the stabiliser of $\xi$ and exchanging orders of summation. Then, in order to find the constants it suffices to choose $\xi=\eta$ and integrate both sides over $\S^n$.
\end{proof}

\begin{lem}[Zonal averaging]
For all $\xi\in\S^n$, $k\in\N_0$ and $i=1,\ldots,d_k(n)$
\begin{equation}\label{eq:zonal-average-Sn}
\int_{\SO(n)} \mathcal Y_k^i(B\xi)\d B = \frac{\mathcal Y_k^i(\xi_0)}{\mathcal C_k^{(n-1)/2}(1)} \mathcal C_k^{(n-1)/2} (\xi_0^\top \xi)
\end{equation}
\end{lem}
\begin{proof}[Sketch of proof]
As the left hand side is after averaging zonal with respect to the north pole $\xi_0$, it is again just a matter of determining constants. For this we set $\xi=\xi_0$.
\end{proof}

\begin{lem}[Funk-Hecke formula]
Let $f:[-1,1]\to\C$ be continuous. Then for all $i=1,\ldots,d_k(n)$
\begin{equation}
   \int_{\S^n} f(\xi^\top\eta) \mathcal Y_k^i(\xi)\d\xi = \mathcal Y_k^i(\eta) 
   \frac{|\S^{n-1}|}{\mathcal C_k^{(n-1)/2}(1)}
    \int_{-1}^1 f(t) \mathcal C_k^{(n-1)/2}(t)(1-t^2)^{n/2-1} \d t.
\end{equation} 
\end{lem}
\begin{proof}[Sketch of proof]
As $f(\xi^\top\eta)$ is zonal with respect to $\eta$, we can average the integrand by replacing $\eta$ with $B\eta$ and integrating over $B$ from the stabiliser of $\eta$. Using the invariance of the integral on $\S^n$ this can be transferred to an average over $\mathcal Y^i_k$ and therefore \eqref{eq:zonal-average-Sn} reduces the integral to an integral over $f(\xi^\top\eta)\mathcal C_k^{(n-1)/2}(\xi^\top\eta)$. It remains to substitute $t=\xi^\top\eta$ and the statement follows.
\end{proof}

\subsubsection{Zonal heat wavelets on $\S^n$}\label{sec:4.3.1} At first we want to give a formula for zonal wavelets. Again we fix $\xi_0$ to be the north pole of $\S^n$ and set $\xi=A\xi_0$. Then the $\SO(n)$-averaged character 
$\trace \mathcal T_k(A)$ has to be zonal w.r.to $\xi_0$ and is therefore given by
\begin{align}
 (\mathbb P_{\S^n}  \trace  \mathcal T_k) (\xi)&=   \frac{\mathcal C_k^{(n-1)/2} (\xi^\top\xi_0)}{\mathcal C^{(n-1)/2}_k(1)}
\trace \mathcal T_k(\I)
 =   \frac{d_k(n)}{\mathcal C_k^{(n-1)/2}(1)}  \mathcal C_k^{(n-1)/2} (\xi^\top\xi_0)\notag\\
 &=\frac{2k+n-1}{n-1}\mathcal C_k^{(n-1)/2} (\xi^\top\xi_0)
\end{align}
based on $d_k(n)=\binom{n+k}{n}-\binom{n+k-2}{n}$ and $C_k^{(n-1)/2}(1)=\binom{n+k-2}{k}$.

Hence, a family of real zonal heat wavelets on $\S^n$ is given by
\begin{equation}
  \psi_\rho(\xi) =\frac1{\sqrt{\alpha(\rho)}} \sum_{k=1}^\infty   \frac{2k+n-1}{n-1}\lambda_k \e^{-  \lambda_k^2\rho/2} C_k^{(n-1)/2}(\xi^\top\xi_0),
\end{equation}
where $-\lambda_k^2$ is the corresponding eigenvalue of the Laplacian on $\SO(n+1)$. 

\begin{rem}
In fact, any other sequence of positive numbers $\lambda_k\to\infty$ will give a suitable wavelet transform with similar properties. In this situation we would replace the heat semigroup with an arbitrary semigroup generated by a non-negative bi-invariant operator (given by the sequence of eigenvalues $\lambda_k^2$ on its eigenspaces $\spann\{\mathcal T_k^{ij} : 1\le i,j \le d_k(n)\})$.
\end{rem}

\subsubsection{Nonzonal wavelets on $\S^n$} To shorten notation let $\H=\SO(n)$. We first calculate the projection $\hat\H(k)$, i.e.,
\begin{align}
\hat\H(k) &= \int_{\SO(n)} \mathcal T_k(B)\d B = \left(\int_{\S^n} \int_{\SO(n)} \mathcal Y^i_k(B^{-1}\xi)\d B\, \mathcal Y^j_k(\xi)\d\xi\right)_{ij}\notag\\
&= \left( \frac{\mathcal Y^i_k(\xi_0)}{\mathcal C_k^{(n-1)/2}(1)} \int_{\S^n} \mathcal C^{(n-1)/2}_k(\xi_0^\top\xi)\mathcal Y^j_k(\xi)\d\xi\right)_{ij}\notag\\
&= \left(\frac{\mathcal Y_k^i(\xi_0)\mathcal Y_k^j(\xi_0)}{[\mathcal C_k^{(n-1)/2}(1)]^2}|\S^{n-1}| \int_{-1}^1 
[\mathcal C_k^{(n-1)/2} (t)]^2 (1-t^2)^{n/2-1}\d t\right)_{ij}\notag\\
&=\frac{|\S^n|}{d_k(n)} \left( \mathcal Y_k^i(\xi_0)\mathcal Y_k^j(\xi_0) \right)_{ij}
\end{align}
based on \eqref{eq:zonal-average-Sn}, Funk-Hecke formula and the normalisation of Gegenbauer polynomials,
cf.~\cite[11.4(17)]{Erdely},
\begin{equation}
   |\S^{n-1}|\int_{-1}^1 [\mathcal C_k^{(n-1)/2} (t)]^2 (1-t^2)^{n/2-1}\d t
   = \frac{[\mathcal C_k^{(n-1)/2}(1)]^2}{d_k(n)}|\S^n|.
\end{equation}
In order to simplify notation,  we {\em assume} that the basis of spherical harmonics $\mathcal Y_k^j(\xi)$ is chosen in such a way that $\mathcal Y_k^1(\xi_0)=\sqrt{d_k(n)/|\S^n|}$ and $\mathcal Y_k^{i}(\xi_0)=0$ for all $i>1$. Then $\hat\H(k)=\diag(1,0,\ldots,0)$. In particular this means
\begin{equation}
\sqrt{\frac{|\S^n|}{d_k(n)}}\,  \mathcal Y_k^i(\xi) = (\mathbb P_{\S^n}\mathcal T_k^{i1})(\xi) = \int_{\SO(n)} \mathcal T_k^{i1}(AB)\d B
  = \mathcal T_k^{i1} (A) , \qquad \xi=A\xi_0. 
\end{equation}

In order to specify a non-zonal wavelet family we have to specify the corresponding Fourier coefficients $\hat\psi_\rho(k)$ for each class-1 representation $\mathcal T_k$. Due to the considerations in Section~\ref{sec:nonzonal} we know that $\hat\psi_\rho(k)$ has entries only in the first row and these entries satisfy \eqref{eq:adm-cond-3mass2}.
To single out an interesting family of non-zonal wavelets we follow \cite{EBCK} and introduce an weight vector
$w(k)\in \C^{d_k(n)}$, $|w(k)|^2=\sum_i |w_i(k)|^2=1$, and define the Fourier coefficients in terms of $w(k)$ as
\begin{equation}
   \hat\psi_\rho(k) = \frac1{\sqrt{\alpha(\rho)}} \lambda_k \e^{-\lambda_k^2\rho/2}\, e_k\otimes w(k),
\end{equation}
where $e_k=(1,0,\ldots,0)^\top\in\C^{d_k(n)}$, $-\lambda_k^2$ is the sequence of eigenvalues of the Laplacian on $\SO(n+1)$ and $\alpha(\rho)$ a suitable weight function. 
The choice $w(k)=e_k$ corresponds to the zonal wavelets constructed in Section~\ref{sec:4.3.1}.

\begin{figure}
\includegraphics[width=11cm]{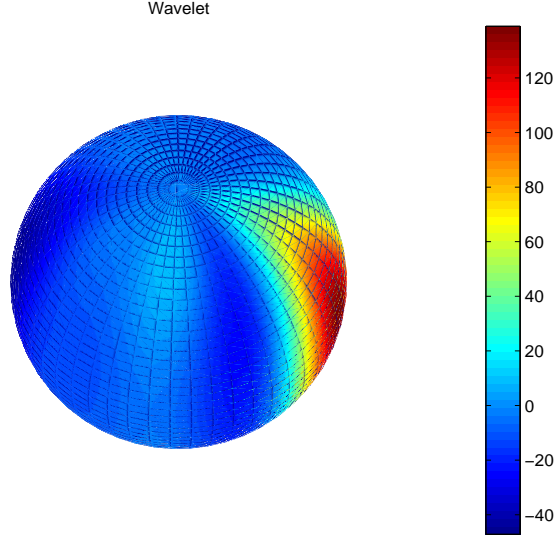}
\caption{An example of a non-zonal heat wavelet on $\S^2$. Scale parameter is $\rho=0.15$, weight-vector 
$w(k) = (2k+1)^{-1/2}(1,\ldots,1)^\top\in \C^{2k+1}$ is equi-distributed and weight function  chosen as $\alpha(\rho)=\rho^{-3}$.}
\label{fig3}
\end{figure}

For $n=2$, $d_k(2)=2k+1$, we get therefore for any choice of weight vectors $w(k)\in \C^{2k+1}$ of unit length
a family of non-zonal wavelets,
\begin{equation}
    \psi_\rho(\xi) = \frac1{\sqrt{\alpha(\rho)}}   \sum_{k=0}^\infty (2k+1) \lambda_k \e^{-\lambda_k^2\rho/2} \sum_{i=-k}^{k} w_{i}(k) \mathcal Y^i_k(\xi),
\end{equation}
where the spherical harmonics on $\S^2$ are denoted in the more conventional way as 
\begin{align}
{\mathcal Y}^0_k(\xi) &=\sqrt{\frac{2k+1}{4\pi}}\, \mathcal C_k^{1/2}(t),\\
{\mathcal Y}^j_k(\xi)&=  
 \sqrt{\frac{(2k+1)\Gamma^2(|j|+\frac12) \Gamma(k-|j|+1)}{2^{2-|j|}\pi^2\Gamma(k+|j|+1)}}
  \mathcal C_{k-|j|}^{|j|+1/2}(t) (1-t^2)^{|j|/2}\e^{\i j \theta},\quad j=-k,\ldots, k,
\end{align}
in the special spherical coordinates $t=\xi_0^\top\xi$ describing the sine of latitude and $\theta\in [-\pi,\pi]$ any convention of longitude on $\S^2$. Figure~\ref{fig3} depicts an example of a non-zonal wavelet on $\S^2$.

The associated non-zonal wavelet transform of a function $f\in L^2(\S^2)$ lives on $\R_+\times \SO(3)$. 
We fix the north pole on $\S^2$ as point with co-ordinates $(0,0,1)^\top$ and introduce Euler angles $\alpha\in[-\pi,\pi]$, $\beta\in[-\pi/2,\pi/2]$, $\gamma\in [-\pi,\pi]$ as coordinates on $\SO(3)$, such that
\begin{align}
A(\alpha,\beta,\gamma)&= 
\begin{pmatrix} \cos\gamma & \sin\gamma & 0\\-\sin\gamma & \cos \gamma & 0 \\ 0&0&1\end{pmatrix}
\begin{pmatrix}1&0&0\\0&\cos\beta&\sin\beta\\0&-\sin\beta&\cos\beta \end{pmatrix}
\begin{pmatrix} \cos\alpha & \sin\alpha & 0\\-\sin\alpha & \cos \alpha & 0 \\ 0&0&1\end{pmatrix}\notag\\
 &=R_\gamma S_\beta R_\alpha.
\end{align} 
Then we can write the wavelet transform as
\begin{align}
\mathscr W f(\rho,\alpha,\beta,\gamma) &= f\bullet \psi_\rho = \int_{\S^2} f(\xi) \overline{\psi_\rho(A^{-1}(\alpha,\beta,\gamma)\xi)}\d\xi\notag\\
&=\int_{\S^2} f( R_\gamma  S_\beta \xi) \overline{\psi_\rho(R_{-\alpha}\xi)}\d\xi.
\end{align}
Hence, $\beta,\gamma$ specify the point under consideration (the preimage of the north pole $(0,0,1)^\top$ under the rotation $R_\gamma S_\beta$), while $\alpha$ specifies directions (related to rotations of the wavelet $\psi_\rho$ around the fixed north pole $(0,0,1)^\top$).

\section{Concluding remarks}

{\bf 1.} 
We defined wavelet transforms for square integrable functions on Lie groups $\G$ and their homogeneous spaces $\mathcal X$. Since heat wavelet families are smooth, all considerations generalise immediately to distributions $u\in \mathcal D'(\G)$ or $u\in\mathcal D'(\mathcal X)$ and give rise to a smooth function $\mathscr Wu(\rho,g)$ decaying exponentially as $\rho\to\infty$. 

Furthermore, the approximation properties of the heat kernel are directly linked to the inversion formula. Hence these are valid in a variety of function spaces, e.g., for any $\phi\in L^p_0(\G)$, $p\in[1,\infty)$,
\begin{equation}
  \phi =  \lim_{t\to0} \int_t^\infty \mathscr W \phi(\rho,\cdot) * \psi_\rho(\cdot) \alpha(\rho)\d\rho\quad\text{in $L^p(\G)$}.
\end{equation}
For arbitrary diffusive approximate identities in the sense of Definition~\ref{df1} only the $L^2$-convergence follows. The transform is still well-defined on $\mathcal D'(\G)$ or $\mathcal D'(\mathcal X)$.

{\bf 2.} 
Regularity properties of functions can be expressed in terms of their wavelet transform. To be more concrete on this, we recall that $\phi\in L^2(\G)$ belongs to the Sobolev space $H^s(\G)$ if the Fourier coefficients satisfy
\begin{equation}
 \|\phi\|_{H^s}^2 =   \sum_{\pi\in\widehat\G} d_\pi \langle\lambda_\pi\rangle^{s} \|\hat\phi(\pi)\|_{HS}^2 < \infty. 
\end{equation}
We say $\phi$ belongs locally in $g\in\G$ to $H^s$, if there exists a cut-off function $\chi\in C^\infty(\G)$ with $g\in\mathrm{supp}\,\chi$ and $\chi\phi\in H^s(\G)$. The latter property can be related to the asymptotic properties of the wavelet transform $\mathscr W\phi(\rho,g)$ evaluated in $g$ as $\rho\to0$. Non-zonal wavelet transforms on homogeneous spaces allow the same consideration micro-localised to directions. 

{\bf 3.} It seems interesting to ask for discretisations of these wavelet transforms. As $\mathscr W\phi(\rho,g)$ is smooth in $\rho>0$ and $g$, we can evaluate it in discrete points. If we choose them in an appropriate way the inversion formul\ae{} imply that the corresponding wavelets form a frame in $L^2(\G)$. We refer to a forthcoming publication concerning details of this.

{\bf 4.} Spherical wavelet transforms are of interest for applications in geophysics, crystallography and medical imaging; see, e.g., \cite{Freeden2} or \cite{Schaeben} and references cited therein. A particularly interesting situation arises for homogeneous spaces $\SO(3) / \Gamma$ for diskrete rotation groups $\Gamma<\SO(3)$ representing crystal symmetries.

\end{document}